\documentclass[10pt, a4paper, reqno]{amsart}

\usepackage[english]{babel}

\usepackage[T1]{fontenc}
\usepackage{lmodern}

\usepackage{amsmath}
\usepackage{amssymb}
\usepackage{amsthm}

\usepackage{url}

\usepackage{hyperref} 

\usepackage{enumitem}
\usepackage[normalem]{ulem}

\newcommand{\lexp}[2]{\null^{#2} \mkern-2mu #1}
\newcommand{\lexpp}[2]{\null^{#2} \mkern-2mu (#1)}

\renewcommand{\epsilon}{\varepsilon}
\renewcommand{\theta}[0]{\vartheta}
\renewcommand{\phi}[0]{\varphi}

\newcommand{\Span}[1]{\left\langle\, #1 \,\right\rangle}

\newcommand{\Set}[1]{\left\{ #1 \right\}}

\DeclareMathOperator{\End}{End}
\DeclareMathOperator{\Aut}{Aut}

\DeclareMathOperator{\Hol}{Hol}
\DeclareMathOperator{\Perm}{Perm}

\newtheorem{dummy}{Dummy}
\numberwithin{dummy}{section}
\numberwithin{figure}{section}

\newtheorem{theorem}[dummy]{Theorem}

\newtheorem{lemma}[dummy]{Lemma}

\newtheorem{proposition}[dummy]{Proposition}

\theoremstyle{definition}
\newtheorem{definition}[dummy]{Definition}

\newtheorem{example}[dummy]{Example}

\theoremstyle{remark}

\newtheorem{remark}[dummy]{Remark}

\makeatletter
\def\imod#1{\allowbreak\mkern10mu({\operator@font mod}\,\,#1)}
\makeatother

\numberwithin{equation}{section}

\usepackage{amsthm, mathtools, amssymb, amsfonts, stmaryrd, latexsym, mathrsfs, todonotes, eufrak, color}

\DeclareMathOperator{\id}{id}

\begin{document}


\date{4 April 2021, 12:12 CEST --- Version 5.00%
}

\author{A. Caranti}

\address[A.~Caranti]%
 {Dipartimento di Matematica\\
  Universit\`a degli Studi di Trento\\
  via Sommarive 14\\
  I-38123 Trento\\
  Italy} 

\email{andrea.caranti@unitn.it} 

\urladdr{https://caranti.maths.unitn.it/}

\author{L. Stefanello}

\address[L.~Stefanello]
        {Dipartimento di Matematica\\
          Universit\`a di Pisa\\
          Largo Bruno Pontecorvo, 5\\
          56127 Pisa\\
          Italy}
\email{lorenzo.stefanello@phd.unipi.it}

\urladdr{https://people.dm.unipi.it/stefanello}
\subjclass[2010]{20B05 20D45 12F10 20N99 16T05 16T25}

\keywords{holomorph, endomorphisms, regular subgroups, skew
          braces, normalising graphs, Yang-Baxter equation, Hopf Galois structures} 

\thanks{The first author is members of INdAM---GNSAGA. The first author
gratefully acknowledges support from the Department of Mathematics of
the University of Trento.}

\title[From endomorphisms to bi-skew braces, etc.]%
      {From endomorphisms to bi-skew braces,\\
        regular subgroups, the Yang--Baxter equation,\\
        and Hopf--Galois structures}

\date{30 June 2021, 11:30 CEST --- Version 6.11}

\begin{abstract}
  The interplay between set-theoretic solutions of the Yang--Baxter
  equation of Mathematical Physics, skew braces, regular subgroups,
  and Hopf--Galois structures has spawned a considerable body of
  literature in recent years.

  In  a recent  paper, Alan  Koch  generalised  a construction  of
  Lindsay N.~Childs,  showing how one can obtain  bi-skew braces $(G,
  \cdot, \circ)$  from an endomorphism  of a group $(G,  \cdot)$ whose
  image is abelian.

  In this  paper, we  characterise the  endomorphisms of  a group
  $(G, \cdot)$ for  which Koch's construction, and a  variation on it,
  yield (bi-)skew braces. We show  how the set-theoretic solutions of the
  Yang--Baxter equation  derived by  Koch's construction carry  over to
  our  more general  situation, and  discuss the  related Hopf--Galois
  structures.
\end{abstract}

\maketitle

\thispagestyle{empty}

\section{Introduction}

The interplay between set-theoretic solutions of the Yang--Baxter
equation of Mathematical Physics, (skew) braces, regular subgroups,
and Hopf--Galois structures has spawned a considerable body of
literature in recent years.

Childs gave in~\cite{Childs-fpf} the  following construction.
Let $G = (G,\cdot)$ be a group, and let $\Perm(G)$ be the group of all
permutations on the underlying set  $G$. Let $\phi$ be an endomorphism
of $G$  which is fixed point  free (meaning that $\lexp{g}{\phi}  = g$
only for $g =  1$; here  where $\lexp{g}{\phi}$ is our notation for
the value of an endomorphism $\phi$ of the group $G$ on an element $g
\in G$),  and abelian (meaning that  the image of $\phi$  is an
abelian  subgroup  of  $G$,  or equivalently  $\lexp{[G,  G]}{\phi}  =
1$). Consider the subset
\begin{equation*}
  N
  =
  \Set{ h \mapsto g \cdot h \cdot (\lexp{g}{\phi})^{-1}
    :
    g \in G
    }
\end{equation*}
of $\Perm(G)$. Then $N$ is a regular subgroup of $\Perm(G)$ which is
normalised by the image $\lambda(G)$ of the left regular representation
of $G$. To such an $N$ is associated a certain group operation $\circ$
on $G$, such that $N$ is isomorphic to $(G, \circ)$, and $(G,
\circ, \cdot)$ is a skew left brace.

In the recent manuscript~\cite{Koch-Abelian} Koch generalised
this construction. Let $\psi$ be an abelian endomorphism of the group
$G$. Then the set
\begin{equation}
  \label{eq:Koch}
  \Set{ h \mapsto g \cdot (\lexp{g}{\psi})^{-1} \cdot h \cdot \lexp{g}{\psi}
    :
    g \in G
    }
\end{equation}
is  a regular  subgroup of  $\Perm(G)$ which  both normalises,  and is
normalised  by, $\lambda(G)$.  The associated  skew brace  $(G, \circ,
\cdot)$ is  thus a  bi-skew brace (\cite{Childs-bi-skew, mybi})  in the
sense    that    $(G,    \cdot,     \circ)$    is    also    a    skew
brace.  In~\cite{Koch-Abelian}  Koch   showed  that  his  construction
encompasses that of Childs.

The goal of this paper is to further extend Koch's construction to the
following situations. Given an endomorphism $\psi$  of $G$ and $\epsilon =
\pm 1$, we consider the  subset
\begin{equation}
  \label{eq:pm1}
  N
  =
  \Set{ h \mapsto
    g \cdot (\lexp{g}{\psi})^{\epsilon} \cdot h \cdot
    (\lexp{g}{\psi})^{-\epsilon} 
    :
    g \in G
    }
\end{equation}
of $\Perm(G)$.  It is immediate to see that $N$ is a regular subset of
$\Perm(G)$ which normalises $\lambda(G)$. We determine the conditions
on  $\psi$ so  that  $N$ is  a  subgroup of  $\Perm(G)$,  and then  is
normalised by $\lambda(G)$.

Our first main result is concerned with the case when $\epsilon = - 1$
in~\eqref{eq:pm1}, which is a direct  extension of the case considered
by Koch.  In the statement, given an endomorphism $\psi$ of $G$ and an
element $g  \in G$,  we  use the
notations $[g, \psi] = g \cdot (\lexp{g}{\psi})^{-1}$, and $[G, \psi]
= \Span{ [g, \psi] : g \in G}$,  the reason being that if $\psi$ is an
automorphism of $G$, then $[g, \psi]$ is indeed the commutator $g
 \psi  g^{-1} 
\psi^{-1}$ in the abstract holomorph of $G$.
\begin{theorem}
  \label{theorem:-1}
  Let $G = (G, \cdot)$ be a group, and let $\psi$ be an endomorphism of $G$.
Consider the subset
  \begin{equation*}
    N
    =
    \Set{ h \mapsto g \cdot
      (\lexp{g}{\psi})^{-1} \cdot
      h \cdot
      \lexp{g}{\psi}
      :
      g \in G
    }
  \end{equation*}
  of $\Perm(G)$. The following are equivalent: \begin{enumerate}[label=(\alph*)]
		\item \label{item:is_a_subgroup}
		$N$ is a subgroup of $\Perm(G)$;
		\item \label{item:is_normalised}
		$N$ is normalised by $\lambda(G)$;
		\item \label{item:is_a_regular_subgroup_plus}
		$N$ is a regular subgroup of $\Perm(G)$ which normalises, and is normalised by, $\lambda(G)$;
		\item \label{item:more_than_Koch}
		$\lexp{[[G,\psi],G]}{\psi}\le Z(G,\cdot)$;
		\item \label{item:biskew}
		$(G, \cdot, \circ)$ is a bi-skew brace, for $g \circ h = g \cdot
      (\lexp{g}{\psi})^{-1} \cdot h \cdot \lexp{g}{\psi}$.
	\end{enumerate}
\end{theorem}
In the situation considered by Koch we have $\lexp{[G, G]}{\psi} = 1$,
so the condition in Theorem~\ref{theorem:-1} is clearly satisfied;  we show
with an example that this  condition is distinct from Koch's condition
$\lexp{[G, G]}{\psi} = 1$.

We then consider the situation when $\epsilon = 1$ in~\eqref{eq:pm1}.
\begin{theorem}
  \label{theorem:+1}
  Let $G=(G,\cdot)$ be a group, and let $\psi$ be an endomorphism of $G$.
  Consider the subset
  \begin{equation*}
    N
    =
    \Set{ h \mapsto g \cdot
      \lexp{g}{\psi} \cdot
      h \cdot
      (\lexp{g}{\psi})^{-1}
      :
      g \in G
    }
  \end{equation*}
  of $\Perm(G)$.

  \begin{enumerate}
  \item
    The following are equivalent:
    \begin{enumerate}
    \item
      \label{item:is_subgroup+}
      $N$ is a subgroup of $\Perm(G)$;
    \item
      \label{item:is_regular_subgroup+}
      $N$ is a regular subgroup of $\Perm(G)$ which normalises
      $\lambda(G)$;
    \item
      \label{item:positive-Koch}
      $\lexp{[\lexp{G}{\psi}, G]}{\psi} \le Z(G)$;
     \item
      \label{item:skew}
      $(G, \cdot, \circ)$ is a skew brace, for
      $g \circ h = g \cdot \lexp{g}{\psi} \cdot h \cdot (\lexp{g}{\psi})^{-1}$.
   \end{enumerate}
  \item
    \label{item:these-are-bi}
    The following are equivalent:
    \begin{enumerate}
    \item
      \label{item:is_bi-subgroup+}
      $N$ is a regular subgroup of $\Perm(G)$ which normalises, and is
      normalised by, $\lambda(G)$;
    \item
      \label{item:positive-Koch-2}
      $\lexp{[G, G]}{\psi} \le Z(G)$;
      \item
      \label{item:biskew-2}
      $(G, \cdot, \circ)$ is a bi-skew brace, for $g \circ h = g \cdot
      \lexp{g}{\psi} \cdot h \cdot (\lexp{g}{\psi})^{-1}$.
  \end{enumerate}
  \end{enumerate}
\end{theorem}
We show with
examples that conditions~\eqref{item:positive-Koch}~and
\eqref{item:positive-Koch-2} of Theorem~\ref{theorem:+1} are distinct, and 
distinct from Koch's condition $\lexp{[G, G]}{\psi} = 1$.

In~\cite{Koch-Abelian} Koch also applied his construction to exhibit
solutions of the Yang--Baxter equation.  We  do the same in the context
of  our Theorems~\ref{theorem:-1}~and  \ref{theorem:+1}.

In~\cite{Koch-Abelian}  Koch derived  results about  the Hopf--Galois
extensions related to his construction.  We are able to obtain similar
results in our setting, studying Hopf algebras associated to regular
subgroups of $\Perm(G)$  which  normalise  $\lambda(G)$.

Section~\ref{sec:preliminaries} deals with some preliminaries. Since
the subject  of skew braces is  fairly recent, some basic  results are
scattered  through  the  literature; we  sum them  up  here  for  the
convenience of  the reader.  

Section~\ref{sec:proofs}
contains the proofs of Theorems~\ref{theorem:-1}~and \ref{theorem:+1},
plus some supplementary material on the structure
of the involved regular
subgroups.
Section~\ref{sec:examples} contains the examples mentioned
above.

Section~\ref{sec:YB}  describes  the  set-theoretic solutions  of  the
Yang--Baxter   equation  associated   to  the   skew braces   of
Theorems~\ref{theorem:-1}~and \ref{theorem:+1}.

In  Section~\ref{sec:HG1}  we  introduce the  Hopf--Galois  structures
associated to  the regular subgroups  of Theorems~\ref{theorem:-1}~and
\ref{theorem:+1},  which we  investigate  in Section~\ref{sec:HG2}  by
considering certain subgroups of the relevant regular subgroups, which
correspond to sub-Hopf algebras.

We are indebted to the referee for their suggestions, which led to a
radical revision of the structure of the paper.

\section{Preliminaries}
\label{sec:preliminaries}

In what follows, if $(G,\cdot)$ is  a group, we write $\End(G, \cdot)$
for the monoid of endomorphisms of $(G, \cdot)$. We denote the action
of  $\psi  \in  \End(G, \cdot)$  on  $g  \in  G$  by a  left  exponent
$\lexp{g}{\psi}$, and  thus compose such  endomorphisms right-to-left.
If $\epsilon$ is an integer, we have  clearly $\lexp{(g^{\epsilon})}{\psi} =
(\lexp{g}{\psi})^{\epsilon}$;      we      thus      write      simply
$\lexp{g}{\psi}^{\epsilon}$ for such an element. Similarly, if $\eta$ an element of $\Perm(G)$, the group of permutations on the underlying set $G$, we write $\lexp{g}{\eta}$ for the image of $g$ under $\eta$.

Let $(G,\cdot)$ be a group.  Denote by $\lambda$ the \emph{left regular representation}: 
\begin{align*}
  \lambda \colon (G, \cdot) &\to \Perm(G)\\
             g &\mapsto (h \mapsto g \cdot h)
\end{align*}

\begin{definition}
The \emph{(permutational) holomorph} of $(G,\cdot)$ is defined as the normaliser of
the image $\lambda(G)$ of $G$ in $\Perm(G)$:
\begin{equation*}
  \Hol(G, \cdot)
  =
  N_{\Perm(G)}(\lambda(G))
  .
\end{equation*}
\end{definition}
It  is  immediate  to  see  that 
\begin{equation*}
	\Hol(G,  \cdot)=\lambda(G) \Aut(G,  \cdot),
\end{equation*}
that is, the holomorph is  the  (inner) semidirect product of $\lambda(G)$ by the
group of automorphisms $\Aut(G, \cdot)$,  the latter being a bona fide
subgroup of $\Perm(G)$.   The group $\Hol(G, \cdot)$  is therefore isomorphic to
the \emph{abstract holomorph} of $(G, \cdot)$, which is the natural (outer)
semidirect product of $(G, \cdot)$ by $\Aut(G, \cdot)$.

Write $1$ for the identity of $(G,\cdot)$.

\begin{definition}
	A subset $N$ of $\Perm(G)$ is \emph{regular} if the map
\begin{align*}
	N&\to G\\
	\eta&\mapsto \lexp{1}{\eta}
\end{align*}
is   bijective.   We   denote   by  $\nu$  the  inverse of this map, so   that
$N=\Set{\nu(g):g\in G}$,  where $\nu(g)$  is the  unique element  of $N$
such that $\lexp{1}{\nu(g)}=g$.
\end{definition}

\begin{example}\label{exa:regularsubset}
	 Let $\gamma :  G \to \Aut(G, \cdot)$ be a function. Then
	 \begin{equation*}
	 	N=\Set{\lambda(g)\gamma(g):g\in G}
	 \end{equation*}
	 is a regular subset of $\Hol(G,\cdot)$, since for every $g\in G$, $\lambda(g)\gamma(g)$ is the unique element taking $1$ to $g$.
\end{example}

Regular subgroups are strictly connected with skew braces.

\begin{definition}
	A \emph{skew (left) brace} is a  triple  $(G,\cdot,\circ)$,  where   $(G,\cdot)$  and
$(G,\circ)$ are groups, and for every $g,h,k\in G$,
\begin{equation*}
  g\circ(h\cdot k) = (g\circ h)\cdot g^{-1}\cdot (g\circ k).
\end{equation*}
(Here $g^{-1}$ denotes the inverse of $g$ with respect to $\cdot$.)
\end{definition}
It is immediate to see that $(G,\cdot)$ and
$(G,\circ)$  share  the   same  identity  $1$ (\cite[Lemma 1.7]{skew}). 

\begin{definition}
  A \emph{bi-skew  brace} is  a triple  $(G,\cdot,\circ)$, where
  both $(G,\cdot,\circ)$ and $(G,\circ,\cdot)$ are skew braces.
\end{definition}

We summarise here the main results relating skew braces and regular subgroups. The first one combines~\cite[Theorem 4.2]{skew} and~\cite[Proposition 4.3]{skew}. The second one can be obtained as the first one, reversing the role of the two operations (see also the proof of~\cite[Proposition 2.1]{Zen19}). Finally, the third one is part of~\cite[Theorem 3.1]{mybi}, translated from right to left.

\begin{theorem}
  Let $(G,\cdot)$ be a group. The following data are equivalent:
  \label{thm:byo}
  \begin{enumerate}
  \item
    an operation $\circ$ on $G$ such that $(G,\cdot,
    \circ)$ is a skew brace;
  \item
    a  regular  subgroup   $N=\Set{\nu(g):g\in G}\le\Perm(G)$ which normalises $\lambda(G)$.
  \end{enumerate}
\end{theorem}

\begin{proposition}
  \label{pro:gp}
  Let   $(G,\cdot)$  be   a  group.   The  following   data  are
  equivalent:
  \begin{enumerate}
  \item  an  operation   $\circ$  such  that  $(G,\circ,
    \cdot)$ is a skew brace;
  \item  a  regular subgroup  $N=\Set{\nu(g):g\in G}\le
    \Perm(G)$ which is normalised by $\lambda(G)$.
  \end{enumerate}
\end{proposition}

\begin{theorem}
  \label{thm:biskew}
  Let $(G, \cdot)$ be a group. The following data are equivalent:
  \begin{enumerate}
  \item
    an  operation $\circ$  on $G$  such that  $(G,\cdot, \circ)$  is a
    bi-skew brace; 
  \item
    a regular  subgroup  $N=\Set{\nu(g):g\in G}\le\Perm(G)$  which
    normalises, and is normalised by, $\lambda(G)$.
    \end{enumerate}
\end{theorem}

\begin{remark}
	In all the previous results, the correspondence is given by 
	\begin{equation*}
		g\circ h=\lexp{h}{\nu(g)}.
	\end{equation*}
	In particular, the map
	\begin{equation*}
		\nu\colon (G,\circ)\to N
	\end{equation*}
	is an isomorphism.
\end{remark}

\section{Proofs of Theorems~\ref{theorem:-1}~and \ref{theorem:+1}}
\label{sec:proofs}

Let $N$ be a subset of $\Perm(G)$ as in~\eqref{eq:pm1}:
\begin{equation*}
  N
  =
  \Set{ h \mapsto
    g \cdot \lexp{g^{\epsilon}}{\psi} \cdot h \cdot
    \lexp{g^{-\epsilon} }{\psi}
    :
    g \in G
    }.
\end{equation*}
Note that every element of $N$ can be written as $\lambda(g)\iota(\lexp{g^{\epsilon}}{\psi})$ for some $g\in G$, where 
\begin{align*}
  \iota :\ &(G, \cdot) \to \Aut(G, \cdot)\\
           &x \mapsto (y \mapsto x \cdot y \cdot x^{-1})
\end{align*}
is the homomorphism mapping $x \in G$ to the conjugation by $x$.
In particular, as in Example~\ref{exa:regularsubset}, $N$ is a regular subset of $\Perm(G)$ which normalises $\lambda(G)$, and we may write $N=\Set{\nu(g):g\in G}$, with
\begin{equation*}
	\nu(g)=\lambda(g)\iota(\lexp{g^{\epsilon}}{\psi}).
\end{equation*}

\subsection{Proof of Theorem~\ref{theorem:-1}}

Let $N=\Set{\nu(g):g\in G}$ with 
\begin{equation*}
	\nu(g)=\lambda(g)\iota(\lexp{g^{-1}}{\psi}).
\end{equation*}
Then $N$ is a regular subset of $\Hol(G,\cdot)$. Let $g,h\in G$. Since 
\begin{equation*}
	\nu(h)^{-1}=(\lambda(h)\iota(\lexp{h^{-1}}{\psi}))^{-1}=\lambda(\lexp{h}{\psi}\cdot h^{-1}\cdot \lexp{h^{-1}}{\psi})\iota(\lexp{h}{\psi})
\end{equation*}
and 
\begin{align*}
	\nu(g)\nu(h)^{-1}&=\lambda(g)\iota(\lexp{g^{-1}}{\psi})\lambda(\lexp{h}{\psi}\cdot h^{-1}\cdot \lexp{h^{-1}}{\psi})\iota(\lexp{h}{\psi})\\
	&=\lambda(g\cdot \lexp{g^{-1}}{\psi}\cdot \lexp{h}{\psi}\cdot h^{-1}\cdot \lexp{h^{-1}}{\psi}\cdot \lexp{g}{\psi})\iota(\lexpp{g^{-1}\cdot h}{\psi})
\end{align*}
we find that $\nu(g)\nu(h)^{-1}\in N$ if and only if \begin{equation*}
		\iota(\lexpp{g\cdot \lexp{g^{-1}}{\psi}\cdot \lexp{h}{\psi}\cdot h^{-1}\cdot \lexp{h^{-1}}{\psi}\cdot \lexp{g}{\psi}}{\psi}^{-1})=\iota(\lexpp{g^{-1}\cdot h}{\psi}),
	\end{equation*}
	that is, \begin{equation*}
 	\lexp{[[h^{-1}\cdot g,\psi],h]}{\psi}=\lexpp{h^{-1}\cdot g\cdot \lexpp{g^{-1}\cdot h}{\psi}\cdot h\cdot \lexpp{h^{-1}\cdot g}{\psi}\cdot g^{-1}\cdot h\cdot h^{-1}}{\psi}\in Z(G,\cdot).
 \end{equation*}
	This shows that~\ref{item:is_a_subgroup} and~\ref{item:more_than_Koch} are equivalent. 
	
	Since
	\begin{equation*}
		\lambda(h)\nu(g)\lambda(h)^{-1}=\lambda(h)\lambda(g)\iota(\lexp{g^{-1}}{\psi})\lambda(h^{-1})=\lambda(h\cdot g\cdot \lexp{g^{-1}}{\psi}\cdot h^{-1}\cdot \lexp{g}{\psi})\iota(\lexp{g^{-1}}{\psi}),
	\end{equation*}
	we find that $\lambda(h)\nu(g)\lambda(h)^{-1}\in N$ if and only if \begin{equation*}
		\iota(\lexpp{h\cdot g\cdot \lexp{g^{-1}}{\psi}\cdot h^{-1}\cdot \lexp{g}{\psi}}{\psi}^{-1})=\iota(\lexp{g^{-1}}{\psi}),
	\end{equation*}
	that is, \begin{equation*}
		\lexp{[[g,\psi],h]}{\psi}=\lexpp{g\cdot \lexp{g^{-1}}{\psi}\cdot h\cdot \lexp{g}{\psi}\cdot g^{-1}\cdot h^{-1}}{\psi}\in Z(G,\cdot).
	\end{equation*}
	This shows that~\ref{item:is_normalised} and~\ref{item:more_than_Koch} are equivalent, and so also equivalent to~\ref{item:is_a_regular_subgroup_plus}.
	
	Finally,~\ref{item:is_a_regular_subgroup_plus} and~\ref{item:biskew} are equivalent by Theorem~\ref{thm:biskew}. 

\subsection{Proof of Theorem~\ref{theorem:+1}}

Let $N=\Set{\nu(g):g\in G}$ with 
\begin{equation*}
	\nu(g)=\lambda(g)\iota(\lexp{g}{\psi}).
\end{equation*}
Then $N$ is a regular subset of $\Hol(G,\cdot)$. Let $g,h\in G$. Since  
\begin{equation*}
	\nu(h)^{-1}=(\lambda(h)\iota(\lexp{h}{\psi}))^{-1}=\lambda(\lexp{h^{-1}}{\psi}\cdot h^{-1}\cdot \lexp{h}{\psi})\iota(\lexp{h^{-1}}{\psi})
\end{equation*}
and 
\begin{align*}
	\nu(g)\nu(h)^{-1}&=\lambda(g)\iota(\lexp{g}{\psi})\lambda(\lexp{h^{-1}}{\psi}\cdot h^{-1}\cdot \lexp{h}{\psi})\iota(\lexp{h^{-1}}{\psi})\\
	&=\lambda(g\cdot \lexp{g}{\psi}\cdot \lexp{h^{-1}}{\psi}\cdot h^{-1}\cdot \lexp{h}{\psi}\cdot \lexp{g^{-1}}{\psi})\iota(\lexpp{g\cdot h^{-1}}{\psi})
\end{align*}
we find that $\nu(g)\nu(h)^{-1}\in N$ if and only if \begin{equation*}
		\iota(\lexpp{g\cdot \lexp{g^{-1}}{\psi}\cdot \lexp{h^{-1}}{\psi}\cdot h^{-1}\cdot \lexp{h}{\psi}\cdot \lexp{g}{\psi}}{\psi})=\iota(\lexpp{g\cdot h^{-1}}{\psi}),
	\end{equation*}
	that is, \begin{equation*}
 	\lexp{[\lexpp{g^{-1}\cdot h^{-1}}{\psi},h^{-1}]}{\psi}\in Z(G,\cdot).
 \end{equation*}
	This shows that~\eqref{item:is_subgroup+} and~\eqref{item:positive-Koch} are equivalent, and so also equivalent to~\eqref{item:is_regular_subgroup+}. Moreover,~\eqref{item:is_regular_subgroup+} and~\eqref{item:skew} are equivalent by Theorem~\ref{thm:byo}.
	
	Since 
	\begin{equation*}
		\lambda(h)\nu(g)\lambda(h)^{-1}=\lambda(h)\lambda(g)\iota(\lexp{g}{\psi})\lambda(h^{-1})=\lambda(h\cdot g\cdot \lexp{g}{\psi}\cdot h^{-1}\cdot \lexp{g^{-1}}{\psi})\iota(\lexp{g}{\psi}),
	\end{equation*}
	we find that $\lambda(h)\nu(g)\lambda(h)^{-1}\in N$ if and only if \begin{equation*}
		\iota(\lexpp{h\cdot g\cdot \lexp{g^{-1}}{\psi}\cdot h^{-1}\cdot \lexp{g}{\psi}}{\psi})=\iota(\lexp{g}{\psi}),
	\end{equation*}
	that is, \begin{equation*}
		\lexpp{h\cdot g\cdot [\lexp{g^{-1}}{\psi},h^{-1}]\cdot h^{-1}\cdot g^{-1}}{\psi}\in Z(G,\cdot).
	\end{equation*}
	If $N$ is a subgroup, then $\lexp{[\lexp{g^{-1}}{\psi},h^{-1}]}{\psi}\in Z(G,\cdot)$, and so we find that $N$ is normalised by $\lambda(G)$ if and only if condition~\eqref{item:positive-Koch-2} holds. 
	As condition~\eqref{item:positive-Koch-2} is stronger than condition~\eqref{item:positive-Koch}, we conclude that~\eqref{item:is_bi-subgroup+} and~\eqref{item:positive-Koch-2} are equivalent.
	
	Finally,~\eqref{item:is_bi-subgroup+} and~\eqref{item:biskew-2} are equivalent by Theorem~\ref{thm:biskew}.

\subsection{Endomorphisms yielding the same subgroup}

In~\cite[Proposition 3.3]{Koch-Abelian} Koch  characterised the pairs
$\psi_1, \psi_2\in\End(G,\cdot)$  of abelian endomorphisms  which give
rise to the  same regular subgroup of $\Perm(G)$. The  very same proof
applies also in our setting.
\begin{proposition}\label{prop:equality}
  Let $(G,\cdot)$ be a group, and let $\psi_1,\psi_2\in\End(G,\cdot)$. Suppose that $\psi_1,\psi_2$ satisfy the condition of Theorem~\ref{theorem:-1} or one of the conditions of Theorem~\ref{theorem:+1}.
 Then the  regular subgroups $N_{1},
  N_{2}$  of  $\Perm(G)$  associated respectively  to  $\psi_1$  and
  $\psi_2$ coincide if and only if, for every $g \in G$,
  \begin{equation*}
    \lexp{g}{\psi_1}\cdot \lexp{g^{-1}}{\psi_2}\in Z(G,\cdot).
  \end{equation*}
\end{proposition}

\subsection{A semidirect decomposition of $N$}

Proposition~\ref{prop:after_Fitting} below is an immediate consequence
of Fitting's  Lemma for  groups, as stated  in \cite[Proof  of Theorem
  4.1]{E-groups} and  \cite[Theorem 4.2]{quasi}, and of  the next
lemma (see~\cite[Section 3]{Koch-two}).
\begin{lemma}
  \label{lem:endocirc}
  In the situation of 
  Theorems~\ref{theorem:-1}~and \ref{theorem:+1}, the endomorphism
  $\psi$ of $(G, \cdot)$ is also an endomorphism of $(G, \circ)$.
\end{lemma}

\begin{proof}
  For every $g, h \in G$, we have
  \begin{equation*}
    \lexp{g}{\psi}\circ \lexp{h}{\psi}= \lexp{g}{\psi}\cdot \lexpp{\lexp{g}{\psi}}{\psi}^{\epsilon}\cdot \lexp{h}{\psi}\cdot \lexpp{\lexp{g}{\psi}}{\psi}^{-\epsilon}
   =\lexpp{g\cdot \lexp{g^{\epsilon}}{\psi}\cdot h\cdot \lexp{g^{-\epsilon}}{\psi}}{\psi} =\lexpp{g\circ h}{\psi}.\qedhere
  \end{equation*}
\end{proof}

Fitting's Lemma for groups now yields the following result.
\begin{proposition}
  \label{prop:after_Fitting}
  Let $(G, \cdot)$ be a group satisfying the descending condition on
  subgroups, and the ascending condition on
  normal subgroups. (In particular, a finite group will do).
 Let $N=\Set{\nu(g): g\in G}$ be one of the regular subgroups arising from
  Theorems~\ref{theorem:-1}~or \ref{theorem:+1}, and let $(G,\cdot, \circ)$
  be the corresponding (bi-)skew brace. Then there is a natural number $n$ such that
  the following hold:
  \begin{enumerate}
  \item
    $J = \ker(\psi^{n})$ is a normal subgroup of both $(G, \cdot)$ and
    $(G, \circ)$;
  \item
    $I =  \lexp{G}{\psi^{n}}$ is a  subgroup of both $(G,  \cdot)$ and
    $(G, \circ)$;
  \item
    $\psi$ restricts to a nilpotent endomorphism on $J$;
  \item
    $\psi$ restricts to an automorphism on $I$;
  \item
    both $(G, \cdot)$ and $(G, \circ)$ are semidirect products of $J$ by $I$. 
  \end{enumerate}
  Applying $\nu$  to $(G, \circ)$, we obtain that
  \begin{enumerate}[resume]
  \item
  $N$ is a semidirect
  product of the normal subgroup $\nu(J)$ by $\nu(I)$.
  \end{enumerate}
\end{proposition}

\section{Examples}
\label{sec:examples}

\subsection{Examples for Theorem~\ref{theorem:-1}}
\begin{example}
	If $G$ is any group and $\psi$ is  the identity  on $G$, then   condition~\ref{item:more_than_Koch}   is
clearly  satisfied, as  $[G,\psi] =  1$. In  this  case, for every $g,h\in G$, 
\begin{equation*}
	g\circ h=g\cdot g^{-1}\cdot h\cdot g=h\cdot g,
\end{equation*}
therefore the  associated skew
brace is  the \emph{almost trivial}  skew brace.
\end{example}

\begin{example}\label{exa:-1}
	Let $S$  be a group  of nilpotence  class two, and let $G =  S \times
S$. Let $\psi : G \to G$ be the projection on the second factor:
\begin{equation*}
  \lexp{(a, b)}{\psi} = (1, b).
\end{equation*}
The  endomorphism $\psi$  is  clearly  not abelian,  as  its image  is
isomorphic to $S$.  For $g = (a, b) \in G$ we have
\begin{equation*}
  [g,\psi] = (a,b) \cdot (1, b^{-1}) = (a, 1),
\end{equation*}
so that
\begin{equation*}
  \lexp{[[G,\psi],G]}{\psi}= 1\le Z(G).
\end{equation*}
Thus $\psi$ satisfies condition~\ref{item:more_than_Koch} of
Theorem~\ref{theorem:-1}, but not Koch's condition $\lexp{[G,
    G]}{\psi} = 1$. 

\end{example}

We now exhibit  an example of a group $G$ with an endomorphism $\psi$ such that $1  \ne \lexp{[[G,\psi],G]}{\psi} \le
Z(G)$, so that in particular $\lexp{[G,  G]}{\psi} \ne 1$.  

\begin{example}
	Let $S$ be
a group of nilpotence class $3$, and let $T = S / [[S, S], S]$, so that $T$
has nilpotence class $2$. Write $\pi  : S \twoheadrightarrow T$ for the
natural epimorphism, and  consider $G = S \times S  \times T$.  Define
$\psi\in \End(G)$ by
\begin{equation*}
  \lexp{(a, b, c)}{\psi} = (a, a, \pi(b)).
\end{equation*}
Again, $\psi$ is non-abelian, as $G$  projects onto $S$.  For $g = (a,
b, c) \in G$ we have
\begin{equation*}
  [g,\psi]
  =
  (a, b, c) \cdot (a^{-1}, a^{-1}, \pi(b^{-1}))
  =
  (1, b \cdot  a^{-1}, c\cdot \pi(b^{-1}))),
\end{equation*}
so that
\begin{equation*}
  [G,\psi] = 1\times S\times T.
\end{equation*}
It follows that   
\begin{align*}
  \lexp{[[G, \psi], G]}{\psi}
  &=
  \lexp{(1 \times [S,S] \times [T,T])}{\psi}
  \\&=
  1 \times  1\times  \pi([S,S])
  \\&=
  1 \times 1 \times  [T, T]
\end{align*}
is a non-trivial subgroup of $Z(G)$.
\end{example}

\subsection{Examples for Theorem~\ref{theorem:+1}}
\begin{example}
	If $G$ is any group and $\psi$ is the identity on $G$, then both
conditions~\eqref{item:positive-Koch}~and \eqref{item:positive-Koch-2}
become $[G, G] \le Z(G)$, that is, they hold when $G$ has nilpotence
class at most two. Here we have, for every $g,h\in G$,
\begin{equation*}
  g \circ h
  =
  g \cdot g\cdot h\cdot g^{-1}
  =
  g \cdot [g,h]\cdot h\cdot g\cdot g^{-1}
  =
  g \cdot [g, h] \cdot h 
  =
  g \cdot h \cdot [g, h].
\end{equation*}
\end{example}

\begin{example}
	If $G$ is a group of nilpotence class two, then any endomorphism
$\psi$ of $G$ satisfies condition~\eqref{item:positive-Koch-2} of
Theorem~\ref{theorem:+1}. If $\psi$ is an automorphism, then it does
not satisfy Koch's condition $\lexp{[G, G]}{\psi} = 1$.
\end{example}

\begin{example}
	Let $S$  be a group  of nilpotence
class $3$. Let $T  = S/[[S, S], S]$, so that  $T$ has nilpotence class
$2$, and let $\pi : S \twoheadrightarrow T$ be the natural projection.
Let $G = S \times T$, and define $\psi\in \End(G)$ by
\begin{equation*}
  \lexp{(x, y)}{\psi} = (1, \pi(x)).
\end{equation*}
Then 
\begin{equation*}
  \lexp{[G, G]}{\psi}
  =
  \lexp{([S, S] \times [T, T])}{\psi}
  =
  1 \times [T, T],
\end{equation*}
is a non-trivial subgroup of $Z(G)$. Thus $G$ satisfies
condition~\eqref{item:positive-Koch-2} of Theorem~\ref{theorem:+1},
but not Koch's condition $\lexp{[G, G]}{\psi} = 1$.
\end{example}

\begin{example}\label{exa:+1}
	Let $S$  be a group of  nilpotence class $3$,
  and  $T  =  S  /[[S,  S],S]$,  $U   =  S  /[S,  S]$,  so  that  $T$,
  respectively\ $U$ have nilpotence  class $2$, respectively $1$.  Let
  $\pi : S \twoheadrightarrow T$ and $\sigma : T \twoheadrightarrow U$
  be the  natural projections.  Take  $G = S  \times T \times  U$, and
  define $\psi \in \End(G)$ by
  \begin{equation*}
    \lexp{(a, b, c)}{\psi}
    =
    (1, \pi(a), \sigma(b)).
  \end{equation*}
  Now
  \begin{equation*}
    \lexp{[G, G]}{\psi}
    =
    \lexp{([S, S] \times [T, T] \times [U, U])}{\psi}
    =
    1 \times [T, T] \times [U, U]=1\times [T,T]\times 1
  \end{equation*}
  is contained in $Z(G)$ and  non-trivial, as $T$ has nilpotence class
  $2$. In particular, $\psi$ is not abelian. However we have
  \begin{equation*}
    \lexp{G}{\psi} = 1 \times T \times U,
  \end{equation*}
  so that
  \begin{equation*}
    \lexp{[\lexp{G}{\psi}, G]}{\psi}
    =
  \lexp{([1, S] \times [T, T] \times [U, U])}{\psi}
  =
  1 \times 1 \times[U, U]
  =
  1.
\end{equation*}
\end{example}

We now construct an example where condition~\eqref{item:positive-Koch} of
Theorem~\ref{theorem:+1} holds,
but not condition~\eqref{item:positive-Koch-2}.

\begin{example}
	Let $S$ be a group of nilpotence class  $4$, let $T = S /[[[S, S], S],
  S]$, and $U =  S /[[S, S], S]$, so that  $T$, respectively\ $U$ have
nilpotence   class   $3$,   respectively\    $2$.   Let   $\pi   :   S
\twoheadrightarrow T$  and $\sigma  : T  \twoheadrightarrow U$  be the
natural projections.  Take $G = S \times T \times U$, and define $\psi
\in \End(G)$ by
\begin{equation*}
  \lexp{(a, b, c)}{\psi}
  =
  (1, \pi(a), \sigma(b)).
\end{equation*}
Now
\begin{equation*}
  \lexp{[G, G]}{\psi}
  =
  \lexp{([S, S] \times [T, T] \times [U, U])}{\psi}
  =
  1 \times [T, T] \times [U, U]
\end{equation*}
is  not  contained  in
$Z(G)$,    as $T$ has nilpotence class $3$. However we have
\begin{equation*}
  \lexp{G}{\psi} = 1 \times T \times U,
\end{equation*}
so that
\begin{equation*}
  \lexp{[\lexp{G}{\psi}, G]}{\psi}
  =
  \lexp{([1, S] \times [T, T] \times [U, U])}{\psi}
  =
  1 \times 1 \times[U, U]
  \le
  Z(G).
  \end{equation*}
\end{example}

\section{Set-theoretic solutions of the Yang--Baxter equation}
\label{sec:YB}

We recall that a \emph{set-theoretic solution of the Yang--Baxter equation}, defined in~\cite{Dri92}, is
a couple  $(X,r)$, where $X\ne\emptyset$  is a set  and
\begin{align*}
  r\colon     X\times     X&\to      X\times     X\\
                      (x,y)&\mapsto  (\sigma_x(y),\tau_y(x))
 \end{align*}
is                  a                   bijective                  map
satisfying                 \begin{equation*}(r\times\id_X)(\id_X\times
  r)(r\times\id_X)=(\id_X\times           r)(r\times\id_X)(\id_X\times
  r).\end{equation*}  
  We say that $(X,r)$ is
  \begin{itemize}
  	\item \emph{non-degenerate} if, for  every $x\in  X$, $\sigma_x$ and $\tau_x$  are
bijective;
	\item \emph{involutive} if $r^2=\id_{X\times X}$. 
  \end{itemize}
In the sequel, we say that
$(X,r)$ is a solution if $(X,r)$
is  a   set-theoretic  non-degenerate  solution  of   the  Yang--Baxter
equation.

In~\cite{braces} Rump found  a connection  between (left)  braces
$(G, \cdot, \circ)$ and involutive solutions. (Recall that a brace can
be  defined  a posteriori  as  a  skew  brace  where $(G,  \cdot)$  is
abelian.)  This  was generalised  in~\cite[Theorem 3.1]{skew} as  follows.  (Recall
that for $g\in G$, the inverse with respect to $\circ$ is denoted by $g^{-1}$, while we write $\overline{g}$ for the inverse with respect to $\circ$.)
 
\begin{theorem}
  Let  $(G,\cdot,\circ)$  be  a  skew brace.  Then  $(G,r)$  is a  solution,
  where
  \begin{align*}
    r \colon G \times G & \to G \times G\\
    (g, h) & \mapsto(g^{-1} \cdot   (g\circ  h),
    \overline{g^{-1}\cdot (g\circ h)} \circ g \circ h).
  \end{align*}
  The solution $(G,r)$ is involutive if and only if $(G,\cdot)$ is abelian.
\end{theorem}

If   $\mathcal{G}=(G,\cdot,\circ)$  is   a   skew   brace,  then   its
\emph{opposite  skew  brace}  (see~\cite{Rum19,KT20a})  is
$\mathcal{G}'=(G,\cdot',\circ)$, where
\begin{equation*}
  g\cdot' h
  =
  h\cdot g.
\end{equation*}
It  follows  that we  get  two  solutions  from  a single  skew  brace
$\mathcal{G}=(G,\cdot,\circ)$,   namely  $(\mathcal{G},r_\mathcal{G})$
and $(\mathcal{G}',r_{\mathcal{G}'})$.  These solutions are nicely related: $r_{\mathcal{G}'}$ is a two-sided inverse of $r_{\mathcal{G}}$, and $r_{\mathcal{G}'}=r_{\mathcal{G}}$ if and only if $(G,\cdot)$ is abelian (see, for instance,~\cite[Theorem 4.1]{KT20a}).

In particular, if $\mathcal{G}=(G,\cdot,\circ)$ is a bi-skew brace, then also $\mathcal{G}_1=(G,\circ,\cdot)$ is a skew brace, and we get four solutions. 

We may summarise all of this as follows.

\begin{proposition}\label{pro:solutions}
	Let $(G,\cdot,\circ)$ be a skew brace. Then we get two solutions:
	\begin{align*}
		r_\mathcal{G}(g,h)&=(g^{-1}\cdot (g\circ h),\overline{g^{-1}\cdot (g\circ h)} \circ g \circ h);\\
		r_{\mathcal{G}'}(g,h)&=((g\circ h)\cdot g^{-1},\overline{(g\circ h)\cdot g^{-1}} \circ g \circ h).
	\end{align*}
	These solutions are one the inverse of the other and coincide if and only if $(G,\cdot)$ is abelian.
	
	If in addition $(G,\cdot,\circ)$ is a bi-skew brace, then we get other two solutions:
	\begin{align*}
		r_{\mathcal{G}_1}(g,h)&=(\overline{g}\circ (g\cdot h),(\overline{g}\circ (g\cdot h))^{-1}\cdot g\cdot h);\\
		r_{\mathcal{G}_1'}(g,h)&=((g\cdot h)\circ \overline{g},((g\cdot h)\circ \overline{g})^{-1}\cdot g\cdot h)
	\end{align*}
	 These solutions are one the inverse of the other and coincide if and only if $(G,\circ)$ is abelian.
\end{proposition}

\subsection{The case $\epsilon=-1$}

This part extends Koch's work; our  solutions coincide
with his in the particular case when the endomorphism $\psi$ is abelian.

\begin{theorem}\label{thm:sol-1}
  Let  $(G,\cdot)$   be  a   group, and let $\psi\in\End(G,\cdot)$. If $\lexp{[[G,\psi],G]}{\psi}\le   Z(G,\cdot)$, then  we   get  four
  solutions:
  \begin{align*}
    r_\mathcal{G}(g,h)&=(\lexp{g^{-1}}{\psi}\cdot h\cdot \lexp{g}{\psi},\lexpp{g^{-1}\cdot h}{\psi}\cdot h^{-1}\cdot \lexp{g}{\psi}\cdot g\cdot  \lexp{g^{-1}}{\psi}\cdot h\cdot \lexpp{h^{-1}\cdot g}{\psi}); 
    \\
    r_{\mathcal{G}'}(g,h)&=(g\cdot\lexp{g^{-1}}{\psi}\cdot h\cdot \lexp{g}{\psi}\cdot g^{-1},\lexp{h}{\psi}\cdot g\cdot\lexp{h^{-1}}{\psi});
    \\
    r_{\mathcal{G}_1}(g,h)&=(\lexp{g}{\psi}\cdot h\cdot 	\lexp{g^{-1}}{\psi},\lexp{g}{\psi}\cdot h^{-1}\cdot \lexp{g^{-1}}{\psi}\cdot
    g\cdot h);
    \\
    r_{\mathcal{G}_1'}(g,h)&=(g\cdot h\cdot \lexp{h^{-1}}{\psi}\cdot g^{-1}\cdot \lexp{h}{\psi},\lexp{h^{-1}}{\psi}\cdot g\cdot \lexp{h}{\psi}).
  \end{align*}
  
  The   solutions  $(G,r_{\mathcal{G}})$   and  $(G,r_{\mathcal{G}'})$
  are one the inverse of the other and
  coincide if and only if $(G,\cdot)$ is abelian.

  The  solutions $(G,r_{\mathcal{G}_1})$  and $(G,r_{\mathcal{G}_1'})$
  are one the inverse of the other and
  coincide     if    and     only     if,     for    every     $g,h\in
  G$,
  \begin{equation*}
    g\cdot    \lexp{g^{-1}}{\psi}\cdot    h\cdot
    \lexp{g}{\psi}=h\cdot        \lexp{h^{-1}}{\psi}\cdot       g\cdot
    \lexp{h}{\psi}.
  \end{equation*}
\end{theorem}

\begin{proof}
By Theorem~\ref{theorem:-1}, $(G,\cdot,\circ)$ is a bi-skew brace, where for every $g,h\in G$, \begin{equation*}
	g\circ h=g\cdot \lexp{g^{-1}}{\psi}\cdot h\cdot \lexp{g}{\psi}.
\end{equation*}
We just need to apply Proposition~\ref{pro:solutions} in this setting. The key observation is the following: if $g,h\in G$, then 
\begin{align*}
    \lexp{(\lexp{g}{\psi}\cdot
      h\cdot\lexp{g^{-1}}{\psi})}{\psi}&=          \lexp{g}{\psi}\cdot
    \lexp{[g^{-1}\cdot\lexp{g}{\psi},h]}{\psi}\cdot      \lexp{(h\cdot
      g^{-1})}{\psi}\\
      &=\lexp{g}{\psi}\cdot
    \lexp{[[g^{-1},\psi],h]}{\psi}\cdot      \lexp{(h\cdot
      g^{-1})}{\psi},
  \end{align*}
  and since $\lexp{[[g^{-1},\psi],h]}{\psi}\in Z(G,\cdot)$, we find that for every $k\in G$, 
  \begin{equation}\label{eq:key-1}
  	\lexpp{\lexp{g}{\psi}\cdot
      h\cdot\lexp{g^{-1}}{\psi}}{\psi}\cdot k\cdot \lexpp{\lexp{g}{\psi}\cdot
      h^{-1}\cdot\lexp{g^{-1}}{\psi}}{\psi}=\lexpp{g\cdot
      h\cdot g^{-1}}{\psi}\cdot k\cdot \lexpp{g\cdot
      h^{-1}\cdot g^{-1}}{\psi}. \tag{$\ast$}
  \end{equation}
  In this proof, every time we employ~\eqref{eq:key-1}, the symbol $\overset{\ast}{=}$ is used. Recall also that $\overline{g}=\lexp{g}{\psi}\cdot g^{-1}\cdot \lexp{g^{-1}}{\psi}$.
  
  The        first         solution        is
  \begin{align*}
    r_\mathcal{G}(g,h)&=(g^{-1}\cdot (g\circ h),\overline{g^{-1}\cdot (g\circ h)} \circ g \circ h)\\
    &=(\lexp{g^{-1}}{\psi}\cdot h\cdot \lexp{g}{\psi},(\overline{\lexp{g^{-1}}{\psi}\cdot h\cdot \lexp{g}{\psi}}) \circ (g \cdot\lexp{g^{-1}}{\psi}\cdot  h\cdot\lexp{g}{\psi})).
  \end{align*}
Here 
\begin{align*}
	\overline{\lexp{g^{-1}}{\psi}\cdot h\cdot \lexp{g}{\psi}}&=\lexpp{\lexp{g^{-1}}{\psi}\cdot h\cdot \lexp{g}{\psi}}{\psi}\cdot \lexp{g^{-1}}{\psi}\cdot h^{-1}\cdot \lexp{g}{\psi}\cdot \lexpp{\lexp{g^{-1}}{\psi}\cdot h^{-1}\cdot \lexp{g}{\psi}}{\psi}\\
	&\overset{\ast}{=}\lexpp{g^{-1}\cdot h\cdot g}{\psi}\cdot \lexp{g^{-1}}{\psi}\cdot h^{-1}\cdot \lexp{g}{\psi}\cdot \lexpp{g^{-1}\cdot h^{-1}\cdot g}{\psi}\\
	&=\lexpp{g^{-1}\cdot h}{\psi}\cdot h^{-1}\cdot \lexpp{h^{-1}\cdot g}{\psi},
\end{align*}
and so
\begin{align*}
	&(\overline{\lexp{g^{-1}}{\psi}\cdot h\cdot \lexp{g}{\psi}}) \circ (g \cdot\lexp{g^{-1}}{\psi}\cdot  h\cdot\lexp{g}{\psi})=\lexpp{g^{-1}\cdot h}{\psi}\cdot h^{-1}\cdot \lexpp{h^{-1}\cdot g}{\psi}\circ(g \cdot\lexp{g^{-1}}{\psi}\cdot  h\cdot\lexp{g}{\psi})\\
	&=\lexpp{g^{-1}\cdot h}{\psi}\cdot h^{-1}\cdot \lexpp{h^{-1}\cdot g}{\psi}\cdot \lexpp{\lexpp{g^{-1}\cdot h}{\psi}\cdot h\cdot \lexpp{h^{-1}\cdot g}{\psi}}{\psi}\\
	&\cdot g \cdot\lexp{g^{-1}}{\psi}\cdot  h\cdot\lexp{g}{\psi}\cdot \lexpp{\lexpp{g^{-1}\cdot h}{\psi}\cdot h^{-1}\cdot \lexpp{h^{-1}\cdot g}{\psi}}{\psi}\\
	&\overset{\ast}{=}\lexpp{g^{-1}\cdot h}{\psi}\cdot h^{-1}\cdot \lexpp{h^{-1}\cdot g}{\psi}\cdot \lexpp{g^{-1}\cdot h\cdot g}{\psi}\cdot g \cdot\lexp{g^{-1}}{\psi}\cdot  h\cdot\lexp{g}{\psi}\cdot \lexpp{g^{-1}\cdot h^{-1}\cdot g}{\psi}\\
	&=\lexpp{g^{-1}\cdot h}{\psi}\cdot h^{-1}\cdot \lexp{g}{\psi}\cdot g \cdot\lexp{g^{-1}}{\psi}\cdot  h\cdot \lexpp{h^{-1}\cdot g}{\psi}.
\end{align*}

The second solution is 
\begin{align*}
	r_{\mathcal{G}'}(g,h)&=((g\circ h)\cdot g^{-1},\overline{(g\circ h)\cdot g^{-1}} \circ g \circ h)\\
	&=(g\cdot \lexp{g^{-1}}{\psi}\cdot h\cdot \lexp{g}{\psi}\cdot g^{-1},\overline{(g\cdot \lexp{g^{-1}}{\psi}\cdot h\cdot \lexp{g}{\psi}\cdot g^{-1})} \circ (g \cdot \lexp{g^{-1}}{\psi} \cdot h\cdot \lexp{g}{\psi})).
\end{align*}
Here 
\begin{align*}
	&\overline{g\cdot \lexp{g^{-1}}{\psi}\cdot h\cdot \lexp{g}{\psi}\cdot g^{-1}}\\&=\lexpp{g\cdot \lexp{g^{-1}}{\psi}\cdot h\cdot \lexp{g}{\psi}\cdot g^{-1}}{\psi}\cdot g\cdot \lexp{g^{-1}}{\psi}\cdot h^{-1}\cdot \lexp{g}{\psi}\cdot g^{-1}\cdot \lexpp{g\cdot \lexp{g^{-1}}{\psi}\cdot h^{-1}\cdot \lexp{g}{\psi}\cdot g^{-1}}{\psi}\\
	&=\lexp{g}{\psi}\cdot \lexpp{\lexp{g^{-1}}{\psi}\cdot h\cdot \lexp{g}{\psi}}{\psi}\cdot \lexp{g^{-1}}{\psi}\cdot g\cdot \lexp{g^{-1}}{\psi}\cdot h^{-1}\cdot \lexp{g}{\psi}\cdot g^{-1}\cdot \lexp{g}{\psi}\cdot\lexpp{ \lexp{g^{-1}}{\psi}\cdot h^{-1}\cdot \lexp{g}{\psi}}{\psi}\cdot \lexp{g^{-1}}{\psi}\\
	&\overset{\ast}{=}
	\lexp{g}{\psi}\cdot \lexpp{g^{-1}\cdot h\cdot g}{\psi}\cdot \lexp{g^{-1}}{\psi}\cdot g\cdot \lexp{g^{-1}}{\psi}\cdot h^{-1}\cdot \lexp{g}{\psi}\cdot g^{-1}\cdot \lexp{g}{\psi}\cdot\lexpp{g^{-1}\cdot h^{-1}\cdot g}{\psi}\cdot \lexp{g^{-1}}{\psi}\\
	&=\lexp{h}{\psi}\cdot g\cdot \lexp{g^{-1}}{\psi}\cdot h^{-1}\cdot \lexp{g}{\psi}\cdot g^{-1}\cdot \lexp{h^{-1}}{\psi},
\end{align*}
and so 
\begin{align*}
	&(\overline{\lexp{g^{-1}}{\psi}\cdot h\cdot \lexp{g}{\psi}}) \circ (g \cdot\lexp{g^{-1}}{\psi}\cdot  h\cdot\lexp{g}{\psi})\\
	&=(\lexp{h}{\psi}\cdot g\cdot \lexp{g^{-1}}{\psi}\cdot h^{-1}\cdot \lexp{g}{\psi}\cdot g^{-1}\cdot \lexp{h^{-1}}{\psi}) \circ (g \cdot\lexp{g^{-1}}{\psi}\cdot  h\cdot\lexp{g}{\psi})\\
	&=\lexp{h}{\psi}\cdot g\cdot \lexp{g^{-1}}{\psi}\cdot h^{-1}\cdot \lexp{g}{\psi}\cdot g^{-1}\cdot \lexp{h^{-1}}{\psi}\cdot\lexpp{\lexp{h}{\psi}\cdot g\cdot \lexp{g^{-1}}{\psi}\cdot h\cdot \lexp{g}{\psi}\cdot g^{-1}\cdot \lexp{h^{-1}}{\psi}}{\psi}\\
	&\cdot g \cdot\lexp{g^{-1}}{\psi}\cdot  h\cdot\lexp{g}{\psi}\cdot \lexpp{\lexp{h}{\psi}\cdot g\cdot \lexp{g^{-1}}{\psi}\cdot h^{-1}\cdot \lexp{g}{\psi}\cdot g^{-1}\cdot \lexp{h^{-1}}{\psi}}{\psi}\\
	&=\lexp{h}{\psi}\cdot g\cdot \lexp{g^{-1}}{\psi}\cdot h^{-1}\cdot \lexp{g}{\psi}\cdot g^{-1}\cdot \lexp{h^{-1}}{\psi}\cdot\lexpp{\lexp{h}{\psi}\cdot g}{\psi}\cdot \lexpp{\lexp{g^{-1}}{\psi}\cdot h\cdot \lexp{g}{\psi}}{\psi}\cdot \lexpp{g^{-1}\cdot \lexp{h^{-1}}{\psi}}{\psi}\\
	&\cdot g \cdot\lexp{g^{-1}}{\psi}\cdot  h\cdot\lexp{g}{\psi}\cdot 
	\lexpp{\lexp{h}{\psi}\cdot g}{\psi}\cdot \lexpp{\lexp{g^{-1}}{\psi}\cdot h^{-1}\cdot \lexp{g}{\psi}}{\psi}\cdot \lexpp{g^{-1}\cdot \lexp{h^{-1}}{\psi}}{\psi}\\
	&\overset{\ast}{=}
	\lexp{h}{\psi}\cdot g\cdot \lexp{g^{-1}}{\psi}\cdot h^{-1}\cdot \lexp{g}{\psi}\cdot g^{-1}\cdot \lexp{h^{-1}}{\psi}\cdot\lexpp{\lexp{h}{\psi}\cdot g}{\psi}\cdot \lexpp{g^{-1}\cdot h\cdot g}{\psi}\cdot \lexpp{g^{-1}\cdot \lexp{h^{-1}}{\psi}}{\psi}\\
	&\cdot g \cdot\lexp{g^{-1}}{\psi}\cdot  h\cdot\lexp{g}{\psi}\cdot 
	\lexpp{\lexp{h}{\psi}\cdot g}{\psi}\cdot \lexpp{g^{-1}\cdot h^{-1}\cdot g}{\psi}\cdot \lexpp{g^{-1}\cdot \lexp{h^{-1}}{\psi}}{\psi}\\
	&=\lexp{h}{\psi}\cdot g\cdot \lexp{g^{-1}}{\psi}\cdot h^{-1}\cdot \lexp{g}{\psi}\cdot g^{-1}\cdot \lexp{h^{-1}}{\psi}\cdot\lexpp{\lexp{h}{\psi}}{\psi}\cdot \lexp{h}{\psi}\cdot \lexpp{\lexp{h^{-1}}{\psi}}{\psi}\\
	&\cdot g \cdot\lexp{g^{-1}}{\psi}\cdot  h\cdot\lexp{g}{\psi}\cdot 
	\lexpp{\lexp{h}{\psi}}{\psi}\cdot \lexp{h^{-1}}{\psi}\cdot \lexpp{\lexp{h^{-1}}{\psi}}{\psi}\\
	&\overset{\ast}{=}
	\lexp{h}{\psi}\cdot g\cdot \lexp{g^{-1}}{\psi}\cdot h^{-1}\cdot \lexp{g}{\psi}\cdot g^{-1}\cdot g \cdot\lexp{g^{-1}}{\psi}\cdot  h\cdot\lexp{g}{\psi}\cdot 
	\lexp{h^{-1}}{\psi}\\
	&=\lexp{h}{\psi}\cdot g\cdot \lexp{h^{-1}}{\psi}.
\end{align*}

The third solution is 
\begin{align*}
	r_{\mathcal{G}_1}(g,h)=(\overline{g}\circ (g\cdot h),(\overline{g}\circ (g\cdot h))^{-1}\cdot g\cdot h).
\end{align*}
Here 
\begin{align*}
	\overline{g}\circ (g\cdot h)&=(\lexp{g}{\psi}\cdot g^{-1}\cdot \lexp{g^{-1}}{\psi})\circ (g\cdot h)\\
	&=\lexp{g}{\psi}\cdot g^{-1}\cdot \lexp{g^{-1}}{\psi}\cdot \lexpp{\lexp{g}{\psi}\cdot g\cdot \lexp{g^{-1}}{\psi}}{\psi}\cdot g\cdot h\cdot \lexpp{\lexp{g}{\psi}\cdot g^{-1}\cdot \lexp{g^{-1}}{\psi}}{\psi}\\
	&\overset{\ast}{=}
	\lexp{g}{\psi}\cdot g^{-1}\cdot \lexp{g^{-1}}{\psi}\cdot \lexpp{g\cdot g\cdot g^{-1}}{\psi}\cdot g\cdot h\cdot \lexpp{g\cdot g^{-1}\cdot g^{-1}}{\psi}\\
	&=\lexp{g}{\psi}\cdot h\cdot \lexp{g^{-1}}{\psi},
\end{align*}
and so 
\begin{align*}
	(\overline{g}\circ (g\cdot h))^{-1}\cdot g\cdot h=\lexp{g}{\psi}\cdot h^{-1}\cdot \lexp{g^{-1}}{\psi}\cdot g\cdot h.
\end{align*}

 Finally, the fourth solution is
 \begin{align*}
  r_{\mathcal{G}_1'}(g,h)&=((g\cdot h)\circ \overline{g},((g\cdot h)\circ \overline{g})^{-1}\cdot g\cdot h).
  \end{align*}
  Here \begin{align*}
  	(g\cdot h)\circ \overline{g}&=g\cdot h\cdot \lexpp{h^{-1}\cdot g^{-1}}{\psi}\cdot \lexp{g}{\psi}\cdot g^{-1}\cdot \lexp{g^{-1}}{\psi}\cdot \lexpp{g\cdot h}{\psi}\\
  	&=g\cdot h\cdot \lexp{h^{-1}}{\psi}\cdot g^{-1}\cdot \lexp{h}{\psi},
  \end{align*}
  and so 
  \begin{equation*}
  	((g\cdot h)\circ \overline{g})^{-1}\cdot g\cdot h=\lexp{h^{-1}}{\psi}\cdot g\cdot \lexp{h}{\psi}.\qedhere
  \end{equation*}
\end{proof}

\subsection{The case $\epsilon=1$}

\begin{theorem}
  Let  $(G,\cdot)$ be  a  group, and let $\psi\in\End(G,\cdot)$.  If
  $\lexp{[\lexp{G}{\psi},G]}{\psi}\le  Z(G,\cdot)$, then  we get
  two solutions:
  \begin{align*}
    r_\mathcal{G}(g,h)
    &=
    (\lexp{g}{\psi}\cdot h\cdot \lexp{g^{-1}}{\psi}, \lexpp{h^{-1}\cdot g}{\psi}\cdot h^{-1}\cdot \lexp{g^{-1}}{\psi}\cdot g\cdot \lexp{g}{\psi}\cdot h\cdot \lexpp{g^{-1}\cdot h}{\psi});
    \\
    r_{\mathcal{G}'}(g,h)
    &=
    (g\cdot \lexp{g}{\psi}\cdot h\cdot \lexp{g^{-1}}{\psi}\cdot g^{-1},\lexpp{g\cdot h^{-1}\cdot g^{-1}}{\psi}\cdot g \cdot  \lexpp{g\cdot h\cdot g^{-1}}{\psi}).
  \end{align*}
  These solutions are one the inverse of the other and coincide if and only if $(G,\cdot)$ is abelian.

  If in  addition $\lexp{[G,G]}{\psi}\le Z(G,\cdot)$, we  get other two
  solutions:
  \begin{align*}
    r_{\mathcal{G}_1}(g,h)
    &=
    (\lexp{g^{-1}}{\psi}\cdot h\cdot \lexp{g}{\psi},
    \lexp{g^{-1}}{\psi}\cdot h^{-1}\cdot \lexp{g}{\psi}\cdot
    g\cdot h);
    \\
    r_{\mathcal{G}_1'}(g,h)&=(g\cdot h\cdot
\lexp{h}{\psi}\cdot g^{-1}\cdot \lexp{h^{-1}}{\psi},
    \lexp{h}{\psi}\cdot g\cdot \lexp{h^{-1}}{\psi}).
  \end{align*}
  These solutions are one the inverse of the other and coincide if and only if, for every $g,h\in G$,
  \begin{equation*}
    g\cdot\lexp{g}{\psi}\cdot h\cdot
    \lexp{g^{-1}}{\psi}=h\cdot\lexp{h}{\psi}\cdot g\cdot
    \lexp{h^{-1}}{\psi}.
  \end{equation*}
\end{theorem}

\begin{proof}
This follows from Theorem~\ref{theorem:+1}, Proposition~\ref{pro:solutions}, and reasonings similar to the ones in the proof of Theorem~\ref{thm:sol-1}.
\end{proof}

\section{Hopf--Galois theory}
\label{sec:HG1}
 Let $L/K$ be a finite field extension. A \emph{Hopf--Galois structure} on $L/K$ consists of a cocommutative $K$-Hopf algebra $H$, together with an $H$-action $\ast$ on $L$, such that $L$ is a left $H$-module algebra, and the $K$-linear map \begin{align*}
		j\colon L\otimes_K H&\to \End_K(L)\\
		x\otimes h&\mapsto (y\mapsto x(h\ast y))
	\end{align*}
	is bijective.
In this situation, we say that the $K$-Hopf algebra $H$ gives a Hopf--Galois structure on $L/K$. We refer to~\cite{ChildsBook} for a general treatment about Hopf--Galois theory.

If $L/K$ is a finite Galois  extension with Galois group $G$, then the following result, known as Greither--Pareigis correspondence (\cite[Theorem 3.1]{GP}), allows  us  to make use of  group theory  in
order to find Hopf--Galois structures on $L/K$.

\begin{theorem}[Greither--Pareigis]
  Let $L/K$ be  a finite Galois extension with Galois  group $G$. Then the Hopf--Galois structures on $L/K$ correspond bijectively to
  the regular subgroups of $\Perm(G)$ normalised by $\lambda(G)$.
\end{theorem} 
 
Explicitly, if $N$  is a regular subgroup of  $\Perm(G)$ normalised by
$\lambda(G)$, then the corresponding Hopf--Galois extension is given by
the $K$-Hopf algebra 
\begin{equation*}
	L[N]^G=\Set{x\in L[N]:g\cdot x=x,\ \forall g\in G},
\end{equation*} where $G$ acts on $L$  via Galois action
and on $N$ via conjugation (after the identification $G\leftrightarrow
\lambda(G)$).
The isomorphism  class of $N$  is called  the \emph{type} of  the Hopf--Galois structure given by $L[N]^G$.
\begin{example}
  If $L/K$ is a finite Galois extension with Galois group
   $(G,\cdot)$  and  
  \begin{align*}
    \rho \colon  (G,\cdot) &\to  \Perm(G)\\
    g &\mapsto (h \mapsto h \cdot g^{-1})
  \end{align*}
  is the  \emph{right regular representation}, then  $\rho(G)$ is a
  regular subgroup  of $\Perm(G)$  normalised by $\lambda(G)$,  which
  corresponds  to  the  classical  Galois structure  given  by  $K[G]$
  (\cite[Proposition 6.10]{ChildsBook}).   Now $\lambda(G)$ is
  regular    and    normalised   by   itself. Since
  $\rho(G)=\lambda(G)$  if  and only  if  $(G,\cdot)$  is abelian,  we
  get that a non-abelian Galois extension has  at least two
  non-isomorphic Hopf--Galois structures. Both are clearly of type $(G,\cdot)$. 
\end{example}

The next result allows us to translate
information  about $N$  to information  about $L[N]^G$.  It is  implicitly  contained in
the proof of~\cite[Theorem 5.2]{GP}, and here it is presented in the formulation of~\cite[Proposition 2.2]{CRVcorrespondence}.

\begin{proposition}
  \label{pro:subalgebras}
  Let $L/K$ be a finite Galois extension with Galois group $G$,
   and let $N\le\Perm(G)$ be a  regular subgroup normalised  by  $\lambda(G)$. Then the  $K$-sub-Hopf  algebras  of  $L[N]^G$
  correspond   bijectively  to the  subgroups   of   $N$  normalised   by
  $\lambda(G)$.
\end{proposition}
Explicitly, if $M$ is a subgroup of $N$ normalised by $\lambda(G)$, then $L[M]^G$ is a $K$-sub-Hopf algebra of $L[N]^G$.

Let us now  recall some facts about opposite  subgroups. If $G$
is a finite group and $N=\Set{\nu(g):g\in G}$ is  a regular subgroup
of $\Perm(G)$ normalised by
$\lambda(G)$, consider the centraliser
\begin{equation*}
  N'=C_{\Perm(G)}(N)\
  =
  \Set{n'\in\Perm(G):
    \eta' \eta = \eta \eta', \text{ for all $\eta \in N$} }.
\end{equation*}
In~\cite[Lemma 2.4.2]{GP} Greither and Pareigis proved that
\begin{equation*}
  N'=\Set{\phi(\eta):   \eta\in    N},
\end{equation*}
where
\begin{equation*}
  \lexp{h}{\phi(\eta)}=\lexp{1}{\nu(h)\eta}.
\end{equation*}
The subgroup $N'$ is called  the \emph{opposite subgroup} of $N$. From
this explicit  description, one can  easily get that $N'$  is regular,
and in  the proof of~\cite[Theorem 2.5]{GP}  it is shown that  $N'$ is
normalised by $\lambda(G)$. In particular, $L[N']^G$ gives a Hopf--Galois structure on $L/K$, which is called the \emph{opposite structure}. Clearly, $N=N'$ if and only if $N$ is abelian.

\begin{example}
  For  a finite Galois extension of fields $L/K$ with Galois group $G$,  
  since $\rho(G) = C_{\Perm(G)}(\lambda(G))$,
  the  Hopf--Galois structure  corresponding  to  $\lambda(G)$ is  the
  opposite  structure  of  the  classical Galois  structure  given  by
  $K[G]$. It is called the \emph{canonical nonclassical Hopf--Galois structure}, and it is studied in details in~\cite{Tru16}.
\end{example}

As one  can expect, the opposite  subgroup is strictly related  to the
opposite skew brace (recall the correspondence between regular subgroups and operations of Proposition~\ref{pro:gp}). 
\begin{proposition}\label{pro:oppositesubsb}
  Let $(G,\cdot)$ be a finite group,  and let $N$ be a regular subgroup
  of $\Perm(G)$ normalised by $\lambda(G)$. If $\mathcal{G}=(G,\circ,\cdot)$ is the
  skew       brace       corresponding        to       $N$,       then
  $\mathcal{G'}=(G,\circ',\cdot)$ is the  skew brace corresponding to
  the opposite subgroup $N'$.
\end{proposition}

\begin{proof}
  Recall that
  \begin{equation*}
    N'=\Set{\phi(\eta):  \eta\in  N},
  \end{equation*}
  where
  \begin{equation*}
    \lexp{h}{\phi(\eta)}=\lexp{1}{\nu(h)\eta}.
  \end{equation*}
  If we now set $ \eta=\nu(g)$, we obtain
  \begin{equation*}
    \lexp{h}{\phi(\eta)}=\lexp{1}{\nu(h)\nu(g)}=\lexp{g}{\nu(h)}=h\circ g,
  \end{equation*}
  that is,
  \begin{equation*}\phi(
    n)=\nu_1(g),
  \end{equation*}
  with $\lexp{h}{\nu_1(g)}=h\circ g=g\circ' h$.
\end{proof}

\begin{remark}
	The proof of Proposition~\ref{pro:oppositesubsb} is  similar to  the  
proof  of~\cite[Proposition
  3.4]{KT20a}, but slightly different,  since in~\cite{KT20a} Koch and
Truman build the skew braces  corresponding to $N$ directly on $N$,
and so their result is up to isomorphism.
\end{remark}

\begin{example}
  Let  $(G,\cdot)$   be  a  finite   group,  and  let $\psi \in   \End(G,\cdot)$ such that $\psi$ satisfies one the the following conditions:
  \begin{enumerate}
  	\item $\lexp{[[G,\psi],G]}{\psi}\le Z(G,\cdot)$ (in this case set $\epsilon=-1$);
  	\item $\lexp{[G,G]}{\psi}\le Z(G,\cdot)$ (in this case set $\epsilon=1$).
  \end{enumerate}
  Then 
  \begin{equation*}
    N=\Set{\nu(g):g\in G}
  \end{equation*}
  is a regular subgroup of $\Perm(G)$ normalised by $\lambda(G)$ (Theorems~\ref{theorem:-1} and~\ref{theorem:+1}), where, for every $g,h\in G$,
  \begin{equation*}
  	\lexp{h}{\nu(g)}=g\circ h= g\cdot \lexp{g^\epsilon}{\psi}\cdot h\cdot \lexp{g^{-\epsilon}}{\psi}.
  \end{equation*}
  We deduce by Proposition~\ref{pro:oppositesubsb} that the opposite subgroup is $N'=\Set{\nu'(g):g\in G}$, 
    with 
    \begin{equation*}\lexp{h}{\nu'(g)} = g\circ' h= h \circ g = 
    h \cdot \lexp{h^{\epsilon}}{\psi}
    \cdot g \cdot \lexp{h^{-\epsilon}}{\psi}.
    \end{equation*}
\end{example}

In~\cite{koch2020skew}  Koch   and  Truman,  given  a   finite  group
$(G,\cdot)$  and  a  regular  subgroup  $N=\Set{\nu(g):g\in G}\le\Perm(G)$  normalised  by
$\lambda(G)$, introduced a  set $\Lambda_N=\lambda(G)\cap    N$ of \emph{$\lambda$-points}
and a set $P_N=\rho(G)\cap   N$ of \emph{$\rho$-points}. The following immediate facts hold:
\begin{itemize}
\item $\Lambda_N$ and $P_N$ are subgroups of $N$;
\item   since  $N$,   $\lambda(G)$,   $\rho(G)$   are  normalised   by
  $\lambda(G)$,  also   $\Lambda_N$  and   $P_N$  are   normalised  by
  $\lambda(G)$;
\item   as    $C_{\Perm(G)}(\rho(G))=\lambda(G)$,   the    action   of
  $\lambda(G)$ on $P_N$ via conjugation is trivial;
 \item $\Lambda_N =\Set{\nu(g):\lambda(g)=\nu(g)}$ and $P_N=\Set{\nu(g): \nu(g)=\rho(g^{-1})}$.
\end{itemize}

\section{Special subgroups of $N$ normalised by $\lambda(G)$, and their
  associated sub-Hopf algebras}

\label{sec:HG2}

If $L/K$ is a finite Galois extension with Galois group $(G,\cdot)$
 and $\psi\in \End(G,\cdot)$ satisfies $\lexp{[[G,\psi],G]}{\psi}\le Z(G,\cdot)$ (in this case set $\epsilon=-1$) or $\lexp{[G,G]}{\psi}\le Z(G,\cdot)$ (in this case set $\epsilon=1$), then $N=\Set{\nu(g):g\in G}$, with
 \begin{equation*}
 	\lexp{h}{\nu(g)}=g\cdot \lexp{g^\epsilon}{\psi}\cdot h\cdot \lexp{g^{-\epsilon}}{\psi},
 \end{equation*} is a regular subgroup of $\Perm(G)$ normalised by $\lambda(G)$ (Theorems~\ref{theorem:-1} and~\ref{theorem:+1}), and so $L[N]^G$ gives a Hopf--Galois structure on $L/K$.
 
The classification of the types of Hopf--Galois structures of a finite field
extension, once the isomorphism type of its Galois group is given, has
met great  interest in recent years.   In this line of  reasoning, one
would hope  to be  able to determine,  given an  endomorphism $\psi\in
\End(G,\cdot)$ as above, the  isomorphism  class of  the
corresponding  regular   subgroup  $N\le\Perm(G)$.   As   pointed  out
in~\cite{Koch-Abelian}, it  seems there is  no easy way to  solve this
question,  even   when  $\psi$   is  an  abelian   endomorphism.  Some
information        about         $N$        is         given        by
Proposition~\ref{prop:after_Fitting}: there exists an integer $n$ such that
\begin{equation*}
  N
  =
  \nu(\ker(\psi^n)) \rtimes \nu(\lexp{G}{\psi^n}).
\end{equation*}
 But this is
not enough, in general, to pinpoint the type of $N$. In the following,
we will thus proceed  as in~\cite[Section 6]{Koch-Abelian} locating
certain   special  subgroups   of   $N$,  which   are  normalised   by
$\lambda(G)$,         and         thus         determining,         by
Proposition~\ref{pro:subalgebras}, $K$-sub-Hopf algebras of $L[N]^G$.

Recall that, for the operation $\circ$ determined by $N$, the map
\begin{align*}
  \nu\colon (G,\circ)&\to N
\end{align*}
is an  isomorphism, so  that (normal) subgroups  of $N$  correspond to
(normal) subgroups of  $(G,\circ)$. Recall furthermore that $\psi$  is also an
endomorphism of $(G,\circ)$ (Lemma~\ref{lem:endocirc}).

\subsection{The case $\epsilon=-1$}

This   case    is   a    direct   generalisation   of    Koch's   work
in~\cite{Koch-Abelian}, and the results we  find coincide with his, in
the particular situation when $\lexp{[G,G]}{\psi}=1$.

Let $L/K$ be a finite Galois extension with Galois group $(G,\cdot)$,
 and let $\psi\in  \End(G,\cdot)$ such
that
\begin{equation*}
  \lexp{[[G,\psi],G]}{\psi}\le Z(G,\cdot).
\end{equation*}
Then $N=\Set{\nu(g): g\in G}$,
with  \begin{equation*}
	\lexp{h}{\nu(g)}=g\cdot \lexp{g^{-1}}{\psi}\cdot h\cdot \lexp{g}{\psi},
\end{equation*}
 is   a  regular
subgroup of $\Perm(G)$ normalised by $\lambda(G)$ 
(Theorem~\ref{theorem:-1}), and so $L[N]^G$ gives a Hopf--Galois structure on $L/K$.

Consider $G_0=\ker(\psi)$. It is normal in $(G,\cdot)$ and in $(G,\circ)$,
therefore its image $N_0= \Set{\nu(g_0): g_0\in G_0}$ is normal in $N$.
 Note that for every $h\in G$, $g_0\in G_0$,
\begin{equation*}
  \lexp{h}{\nu(g_0)}
  =
  g_0\cdot\lexp{g_0^{-1}}{\psi}\cdot h\cdot \lexp{g_0}{\psi}
  =
  g_0\cdot h
  =
  \lexp{h}{\lambda(g_0)},
\end{equation*}
that  is, $N_0=\Set{\lambda(g_0):  g_0\in  G_0}=\lambda(G_0)$, and  so
$N_0$ is  also isomorphic to $  (G_0,\cdot)$, and it is  normalised by
$\lambda(G)$.

Consider now the subgroup $G_1=\Set{g\in G: \lexp{g}{\psi}=g}$ of 
$(G,\cdot)$ and $(G,\circ)$ of fixed points under $\psi$.   
In general, this is not a normal subgroup of $(G,\circ)$. Define
$N_1=\Set{\nu(g_1): g_1\in G_1}$. For every $h\in G$, $g_1\in G_1$, we
have
\begin{equation*}\lexp{h}{\nu(g_1)}=g_1\cdot
  g_1^{-1}\cdot                    h\cdot                   g_1=h\cdot
  g_1=\lexp{h}{\rho(g_1^{-1})},
\end{equation*}
that  is,  $N_1=\rho(G_1)$,  and  so   $N_1$  is  also  isomorphic  to
$(G_1,\cdot)$, and it is normalised by $\lambda(G)$.
	
Clearly $G_0\circ G_1=G_0\cdot G_1$.  Since $N_0$ is normal in $N$,  
\begin{equation*}
N_{01}\coloneqq N_0 N_1=\nu(G_0)\nu(G_1)=\nu(G_0\circ G_1)=\nu(G_0\cdot G_1)	
\end{equation*}
 is a subgroup
of $N$.   Note that $N_1$  is normal in $N_{01}$:  as in~\cite[Proposition
  6.3]{Koch-Abelian},    for   every    $g_0\in   G_0$,    $g_1,h_1\in
G_1$, we have
\begin{align*}
  \nu(g_0)\nu(g_1)\nu(h_1)(\nu(g_0)\nu(g_1))^{-1}
  &=
  \lambda(g_0)\rho(g_1^{-1})\rho(h_1^{-1})(\lambda(g_0)\rho(g_1^{-1}))^{-1}
  \\&=
  \lambda(g_0)\rho(g_1^{-1})\rho(h_1^{-1})\rho(g_1)\lambda(g_0^{-1})
  \\&=
  \rho(g_1^{-1}\cdot h_1^{-1}\cdot g_1)\in \rho(G_1)=N_1.
\end{align*}

Since  $G_0\cap G_1$  is trivial,  we  conclude that  $N_{01}$ is  the
direct  product of  $N_0$  and $N_1$,  it  is isomorphic  to
$(G_0,\cdot)\times (G_1,\cdot)$, and it is normalised by $\lambda(G)$.
	
We  can also  use the  $\lambda$-points and  $\rho$-points to  find
other subgroups normalised by $\lambda(G)$: 
\begin{align*}
	\Lambda_N&=\Set{\nu(g):\nu(g)=\lambda(g)}=\Set{\nu(g): \lambda(g)\iota(\lexp{g^{-1}}{\psi})=\lambda(g)}\\
	&=\Set{\nu(g):\iota(\lexp{g^{-1}}{\psi})=\id_G}
	=\Set{\lambda(g):
      \lexp{g}{\psi}\in        Z(G,\cdot)};\\
    P_N&=\Set{\nu(g): \nu(g)=\rho(g^{-1})}= \Set{\nu(g): \lambda(g)\iota(\lexp{g^{-1}}{\psi})=\rho(g^{-1})}\\
    &=\Set{\nu(g):\iota(g\cdot \lexp{g^{-1}}{\psi})=\id_G}
	=\Set{\rho(g): g\cdot\lexp{g^{-1}}{\psi}\in Z(G,\cdot)}.
\end{align*}

Since $\lambda(G)$ acts trivially on  $P_N$ and $N_1\subseteq P_N$,  $\lambda(G)$ acts trivially also on $N_1$.

Summarising, by Proposition~\ref{pro:subalgebras}, we get (up to) five
$K$-sub-Hopf algebras of $L[N]^G$:
\begin{itemize}
\item $L[N_0]^{G}$; 
\item $L[N_1]^{G}=K[N_1]$;
\item $L[N_{01}]^{G}\cong (L[N_0]\otimes_L L[N_1])^G\cong L[N_0]^G\otimes_K K[N_1]$ (this $K$-Hopf algebra isomorphism follows, for example, by Galois descent: see~\cite[2.12]{ChildsBook});
\item $L[\Lambda_N]^{G}$;
\item $L[P_N]^{G}=K[P_N]$.
\end{itemize}

\begin{remark}\label{rk:different-1}
  Some  of these  $K$-sub-Hopf  algebras  may coincide.  For example,  if
  $Z(G,\cdot)=1$, then $N_0=\Lambda_N$ and $N_1=P_N$. But they may also
  be all distinct. Consider Example~\ref{exa:-1}:
  if $S$ is a group of nilpotence class two and $G=S\times S$, we can
  define $\psi : G \to G$ to be the projection on the second factor:
\begin{equation*}
  \lexp{(a, b)}{\psi} = (1, b).
\end{equation*}
Then $1\ne Z(S)\ne S$, $\lexp{[[G,\psi],G]}{\psi}\le Z(G)$, and \begin{itemize}
	\item $N_0=\nu(\Set{(a,b)\in G: b=1})=\nu(S\times 1)$;
	\item $N_1=\nu(\Set{(a,b)\in G: a=1})=\nu(1\times S)$;
	\item $N_{01}=\nu((S\times1) (1\times S))=\nu(G)$;
	\item $\Lambda_N=\nu(\Set{(a,b)\in G: b\in Z(S)} )=\nu(S\times Z(S))$;
	\item $P_N=\nu(\Set{(a,b): a\in Z(S)} )=\nu(Z(S)\times S)$.
\end{itemize}
Since all these subgroups of $N$ are distinct, 
by Proposition~\ref{pro:subalgebras}
they yield distinct $K$-sub-Hopf algebras of $L[N]^G$. 
\end{remark}

Finally, we discuss about situations in which the type of the structure given by $L[N]^G$ is explicit.

If $\psi$ is different from zero and idempotent, then, for every $n\ge
1$,         $\psi^n=\psi$.           In         particular,         by
Proposition~\ref{prop:after_Fitting},
\begin{equation*}
  N
  =
  \nu(\ker(\psi))\rtimes \nu(\lexp{G}{\psi}).
\end{equation*}
As above,  $G_0=\ker(\psi)$ and      $G_1=\Set{g\in       G:
  \lexp{g}{\psi}=g}$. Since $\psi^2=\psi$, it immediately follows that
$G_1=\lexp{G}{\psi}$, that is,
\begin{equation*}
  N
  =
  \nu(G_0)\rtimes \nu(G_1)
  =
  N_0\rtimes N_1.
\end{equation*}
We have seen that this product is actually direct, and $N_0\cong (G_0,\cdot )$, $N_1
\cong (G_1,\cdot)$, so we conclude that 
\begin{equation*}
  N
  \cong 
  (\ker(\psi),\cdot)\times (\lexp{G}{\psi},\cdot).
\end{equation*}

Now suppose that $\psi$ is fixed point free. If $\psi$ is also abelian, then by~\cite{Childs-fpf} and~\cite[Section 4]{Koch-Abelian}, $N\cong(G,\cdot)$. 
Here instead of $\lexp{[G,G]}{\psi}=1$, we may assume the weaker condition $\lexp{[[G,\psi],G]}{\psi}=1$ (see Example~\ref{exa:-1}), and still find that $N\cong(G,\cdot)$. Indeed, since $\psi$ is fixed point free, then 
\begin{align*}
	\alpha\colon G&\to G\\
	g&\mapsto g\cdot \lexp{g^{-1}}{\psi},
\end{align*}
is bijective, and we claim that $\alpha\colon (G,\circ)\to (G,\cdot)$ is an isomorphism. For every $g,h\in G$, we have
  \begin{align*}
    \lexpp{g\circ h}{\alpha}&= \lexpp{g\cdot\lexp{g^{-1}}{\psi}\cdot h\cdot \lexp{g}{\psi}}{\alpha}\\
    &=g\cdot\lexp{g^{-1}}{\psi}\cdot h\cdot \lexp{g}{\psi}\cdot \lexpp{\lexp{g^{-1}}{\psi}\cdot h^{-1}\cdot \lexp{g}{\psi}\cdot g^{-1}}{\psi}\\
    &=g\cdot\lexp{g^{-1}}{\psi}\cdot h\cdot  \lexpp{g\cdot \lexp{g^{-1}}{\psi}\cdot h^{-1}\cdot \lexp{g}{\psi}\cdot g^{-1}\cdot h}{\psi}\cdot \lexp{h^{-1}}{\psi}\\
    &=g\cdot\lexp{g^{-1}}{\psi}\cdot h\cdot  \lexp{[[g,\psi],h^{-1}]}{\psi}\cdot \lexp{h^{-1}}{\psi}=g\cdot\lexp{g^{-1}}{\psi}\cdot h\cdot \lexp{h^{-1}}{\psi}\\
    &=\lexp{g}{\alpha}\cdot \lexp{h}{\alpha}.
  \end{align*}
  Since $\nu\colon (G,\circ)\to N$ is an isomorphism, we derive our assertion. In particular, $N$ and $\lambda(G)$ are isomorphic: 
\begin{equation*}
  \phi\colon N\xrightarrow{\nu^{-1}}  (G,\circ)\xrightarrow{\alpha}
  (G,\cdot)\xrightarrow{\lambda}\lambda(G).
\end{equation*} 
Under   this isomorphism,  an   element  $\nu(g)$   is  sent   to
$\lambda(g \cdot \lexp{g^{-1}}{\psi})$.
An  isomorphism of regular subgroups  of $\Perm(G)$
normalised by $\lambda(G)$ yields  an isomorphism of the corresponding
Hopf  algebras  if  and  only  it  it  is  $G$-equivariant  (see,  for
instance,~\cite[Corollary 2.3]{KKTU19a}), where the $G$-action is  via
conjugation,     after     the    identification     $G\leftrightarrow
\lambda(G)$.

We claim that  $\phi$ yields  an  isomorphism  $L[N]^G\cong
L[\lambda(G)]^G$  as $K$-Hopf  algebras.    
We  need to  check whether,  for
every $g,h\in G$,
\begin{equation}\label{eq:isoha-1}
  \phi(\lambda(h)\nu(g)\lambda(h)^{-1})
  =
  \lambda(h)\phi(\nu(g))\lambda(h^{-1}).
\end{equation}
The  right-hand   side  is  \begin{equation*}   \lambda(h\cdot  g\cdot
  \lexp{g^{-1}}{\psi}\cdot h^{-1}).
\end{equation*}
Since 
$\lambda(h)\nu(g)\lambda(h)^{-1}=\nu(h\cdot   g\cdot   \lexp{g^{-1}}{\psi}\cdot   h^{-1}\cdot
  \lexp{g}{\psi})$
(we have already performed this computation in the proof of Theorem~\ref{theorem:-1}),
the left-hand side is
\begin{align*}
  &\phi(\nu(h\cdot   g\cdot   \lexp{g^{-1}}{\psi}\cdot   h^{-1}\cdot
  \lexp{g}{\psi}))
  \\&=
  \lambda(h\cdot   g\cdot   \lexp{g^{-1}}{\psi}\cdot   h^{-1}\cdot
  \lexp{g}{\psi})\lambda(\lexpp{h\cdot   g\cdot   \lexp{g^{-1}}{\psi}\cdot   h^{-1}\cdot
  \lexp{g}{\psi}}{\psi}^{-1})
  \\&=
  \lambda(h\cdot   g\cdot   \lexp{g^{-1}}{\psi}\cdot   h^{-1}\cdot
  \lexp{g}{\psi})\lambda(\lexpp{\lexp{g^{-1}}{\psi}\cdot h\cdot \lexp{g}{\psi}\cdot g^{-1}\cdot h^{-1}}{\psi})\\
  &=\lambda(h\cdot   g\cdot   \lexp{g^{-1}}{\psi}\cdot   h^{-1})\lambda(\lexpp{g\cdot \lexp{g^{-1}}{\psi}\cdot h\cdot \lexp{g}{\psi}\cdot g^{-1}\cdot h^{-1}}{\psi})\\
  &=\lambda(h\cdot   g\cdot   \lexp{g^{-1}}{\psi}\cdot   h^{-1})\lambda(\lexp{[[g,\psi],h]}{\psi})=\lambda(h\cdot   g\cdot   \lexp{g^{-1}}{\psi}\cdot   h^{-1}),
\end{align*}
and so~\eqref{eq:isoha-1} holds.

\subsection{The case $\epsilon=1$}

Let $L/K$ be a finite Galois extension with Galois group $(G,\cdot)$,
 and let $\psi\in  \End(G,\cdot)$ such
that
\begin{equation*}
  \lexp{[G,G]}{\psi}\le Z(G,\cdot).
\end{equation*}
Then $N=\Set{\nu(g): g\in G}$,
with  \begin{equation*}
	\lexp{h}{\nu(g)}=g\cdot \lexp{g}{\psi}\cdot h\cdot \lexp{g^{-1}}{\psi},
\end{equation*}
 is   a  regular
subgroup of $\Perm(G)$ normalised by $\lambda(G)$ 
(Theorem~\ref{theorem:+1}), and so $L[N]^G$ gives a Hopf--Galois structure on $L/K$.

As  above, $N_0=\Set{\nu(g):g\in  \ker(\psi)}$ equals  $\lambda(G_0)$,
and  it  is   a  normal  subgroup  of  $N$  which   is  normalised  by
$\lambda(G)$, yielding the $K$-sub-Hopf algebra $L[N_0]^{G}$.

However,  in  this  case,  $N_1  =  \Set{\nu(g):  \lexp{g}{\psi}=g}  =
\Set{\lambda(g^2)\rho(g):\lexp{g}{\psi}  =  g}$  is not  (in  general)
normalised  by $\lambda(G)$:  if $g,h\in  G$ with  $\lexp{g}{\psi}=g$,
then
\begin{equation*}
  \lambda(h)\lambda(g^2)\rho(g)\lambda(h^{-1})=\lambda(h\cdot
  g^{2}\cdot h^{-1})\rho(g)
\end{equation*}
belongs to $N_1$ if and only if
\begin{equation*}
  h\cdot g^2\cdot h^{-1}=g^2.
\end{equation*}
This happens, for example, if $(G,\cdot)$ is abelian, but it is false
in general.

We can
find once more explicitly the $\lambda$-points and $\rho$-points:
\begin{align*}
	\Lambda_N&=\Set{\nu(g):\nu(g)=\lambda(g)}=\Set{\nu(g): \lambda(g)\iota(\lexp{g}{\psi})=\lambda(g)}\\
	&=\Set{\nu(g): \iota(\lexp{g}{\psi})=\id_G}=\Set{\lambda(g):
      \lexp{g}{\psi}\in        Z(G,\cdot)};\\
    P_N&=\Set{\nu(g): \nu(g)=\rho(g^{-1})}= \Set{\nu(g): \lambda(g)\iota(\lexp{g}{\psi})=\rho(g^{-1})}\\
    &=\Set{\nu(g):\iota(g\cdot \lexp{g}{\psi})=\id_G}=\Set{\rho(g): g\cdot\lexp{g}{\psi}\in Z(G,\cdot)}.
\end{align*}
	
Since the action of $\lambda(G)$  via conjugation on $P_N$ is trivial,
we  find   (up  to)  three  $K$-sub-Hopf algebras  of
$L[N]^G$:
\begin{itemize}
\item $L[N_0]^{G}$; 
\item $L[\Lambda_N]^{G}$;
\item $L[P_N]^{G}=K[P_N]$.
\end{itemize}

\begin{remark}
	If $Z(G,\cdot)=1$, then $N_0=\Lambda_N$,
	hence  $L[N_0]^{G}=L[\Lambda_N]^{G}$. But these sub-Hopf algebras
	may be all different,
	as the same example of Remark~\ref{rk:different-1} immediately shows. 
\end{remark}

We  conclude  this  subsection  with  a  study  of  fixed  point  free
endomorphisms. Suppose that $\psi$ is fixed point free. We claim that if $\lexp{[\lexp{G}{\psi},G]}{\psi}=1$, then $N\cong(G,\cdot)$ (note that this condition is weaker than the condition $\lexp{[G,G]}{\psi}=1$, as Example~\ref{exa:+1} shows). Indeed, since $\psi$ is fixed point free, the map
  \begin{align*}
    \alpha\colon G
    &
    \to G
    \\
    g
    &\mapsto g\cdot \lexp{g}{\psi}
  \end{align*}
  is bijective. As $N\cong (G,\circ)$ via  $\nu$, it is enough to show
  that
  \begin{equation*}
    \alpha\colon (G,\circ)\to (G,\cdot)
  \end{equation*}
  is a homomorphism.
  For every $g,h\in G$, we have
\begin{align*}
    \lexpp{g\circ h}{\alpha}&= \lexpp{g\cdot\lexp{g}{\psi}\cdot h\cdot \lexp{g^{-1}}{\psi}}{\alpha}\\
    &=g\cdot\lexp{g}{\psi}\cdot h\cdot \lexp{g^{-1}}{\psi}\cdot \lexpp{g\cdot\lexp{g}{\psi}\cdot h\cdot \lexp{g^{-1}}{\psi}}{\psi}\\
    &=g\cdot\lexp{g}{\psi}\cdot h\cdot \lexpp{\lexp{g}{\psi}\cdot h\cdot \lexp{g^{-1}}{\psi}\cdot h^{-1}}{\psi}\cdot \lexp{h}{\psi}\\
    &=g\cdot\lexp{g}{\psi}\cdot h\cdot \lexp{[\lexp{g}{\psi},h]}{\psi}\cdot \lexp{h}{\psi}=g\cdot\lexp{g}{\psi}\cdot h\cdot  \lexp{h}{\psi}\\
    &=\lexp{g}{\alpha}\cdot \lexp{h}{\alpha}.
  \end{align*}
In particular, $N$ and $\lambda(G)$ are isomorphic: 
\begin{equation*}
  \phi\colon N\xrightarrow{\nu^{-1}}  (G,\circ)\xrightarrow{\alpha}
  (G,\cdot)\xrightarrow{\lambda}\lambda(G).
\end{equation*} 
Under   this isomorphism,  an   element  $\nu(g)$   is  sent   to
$\lambda(g \cdot \lexp{g}{\psi})$. Finally, we claim that $\phi$ yields an isomorphism 
 $L[N]^G\cong L[\lambda(G)]^G$  as $K$-Hopf  algebras if and only if $\psi$ is abelian.    
We  need to  check whether, for every $g,h\in G$,
\begin{equation}\label{eq:isoha+1}
  \phi(\lambda(h)\nu(g)\lambda(h)^{-1})
  =
  \lambda(h)\phi(\nu(g))\lambda(h^{-1}).
\end{equation}
The  right-hand   side  is  \begin{equation*}   \lambda(h\cdot  g\cdot
  \lexp{g}{\psi}\cdot h^{-1}).
\end{equation*}
Since 
$
	\lambda(h)\nu(g)\lambda(h)^{-1}=\nu(h\cdot   g\cdot   \lexp{g}{\psi}\cdot   h^{-1}\cdot
  \lexp{g^{-1}}{\psi})$
(we have already performed this computation in the proof of Theorem~\ref{theorem:+1}),
the left-hand side is
\begin{align*}
  &\phi(\nu(h\cdot   g\cdot   \lexp{g}{\psi}\cdot   h^{-1}\cdot
  \lexp{g^{-1}}{\psi} ))
  \\&=
  \lambda(h\cdot g\cdot \lexp{g}{\psi}\cdot
  h^{-1}\cdot     \lexp{g^{-1}}{\psi})\lambda(\lexp{(h\cdot     g\cdot
    \lexp{g}{\psi}\cdot                                    h^{-1}\cdot
    \lexp{g^{-1}}{\psi})}{\psi})
  \\&=
  \lambda(h\cdot       g\cdot
  \lexp{g}{\psi}\cdot          h^{-1})\lambda(\lexpp{g^{-1}\cdot h\cdot g\cdot h^{-1}\cdot h\cdot \lexp{g}{\psi}\cdot h^{-1}\cdot \lexp{g^{-1}}{\psi}}{\psi})\\
  &=\lambda(h\cdot       g\cdot
  \lexp{g}{\psi}\cdot          h^{-1})\lambda(\lexp{[g^{-1},h]}{\psi}\cdot \lexp{[h,\lexp{g}{\psi}]}{\psi})\\
  &=\lambda(h\cdot       g\cdot
  \lexp{g}{\psi}\cdot          h^{-1})\lambda(\lexp{[g^{-1},h]}{\psi}).
\end{align*}
Therefore~\eqref{eq:isoha+1} holds if and only if $\psi$ is abelian.

 

\providecommand{\bysame}{\leavevmode\hbox to3em{\hrulefill}\thinspace}
\providecommand{\MR}{\relax\ifhmode\unskip\space\fi MR }
\providecommand{\MRhref}[2]{%
  \href{http://www.ams.org/mathscinet-getitem?mr=#1}{#2}
}
\providecommand{\href}[2]{#2}

\end{document}